\documentclass[11pt]{amsart}
\usepackage{amssymb}

\usepackage{amscd}

\newtheorem{lemma}{Lemma}[section]

\newcommand{\BN}{{\mathbb{N}}}
\newcommand{\BR}{{\mathbb{R}}}

\newcommand{\BF}{{\mathbb{F}}}

\newcommand{\BP}{{\mathbb{P}}}

\newcommand{\BG}{{\mathbb{G}}}

\newcommand{\BK}{{\mathbb{K}}}

\newcommand{\BSL}{{\mathbb{SL}}}
\newcommand{\OO}{{\mathcal{O}}}
\newcommand{\gD}{\Delta}

\newcommand{\gd}{\delta}
\newcommand{\gb}{\beta}

\newcommand{\gC}{\Gamma}
\newcommand{\gc}{\gamma}
\newcommand{\gs}{\sigma}
\newcommand{\gS}{\Sigma}
\newcommand{\gO}{\Omega}
\newcommand{\gep}{\epsilon}

\newcommand{\ga}{\alpha}
\newcommand{\gt}{\tau}

\newcommand{\SL}{\text{SL}}
\newcommand{\GL}{\text{GL}}

\newtheorem{prop}{Proposition}[section]
\newtheorem{thm}[prop]{Theorem}
\newtheorem{lem}[prop]{Lemma}
\newtheorem{sublem}[prop]{Sublemma}
\newtheorem{cor}[prop]{Corollary}

\theoremstyle{definition}

\newtheorem{rem}[prop]{Remark}

\begin{document}
\author{T. Gelander}
\date{March 17, 2008}
\title{Notes on uniform exponential growth and Tits alternative}
\maketitle

%\begin{abstract}\gD_k

%\end{abstract}

These notes contain results concerning uniform exponential growth which were obtained in collaborations with E. Breuillard and A. Salehi-Golsefidy, mostly during 2005, improving Eskin-Mozes-Oh theorem \cite{EMO}, as well as a uniform uniform version of Tits alternative improving \cite{uti}.

\section{Main statements about Entropy}
Let $\gS$ be a subset of a group $\gC$. The algebraic entropy of $\gS$ is defined to be
$$
 h(\gS):=\lim_{n\to\infty}\frac{1}{n}(\gS\cup 1)^n.
$$
If $h(\gS)>0$, we say that $\gS$ has exponential growth. If $\inf\{h(\gS):\gC=\langle\gS\rangle\}>0$, we say that $\gC$ has uniform exponential growth. We denote by $d^+(\gS)\in\BN\cup\{\infty\}$ the minimal integer $n$, if such exists, for which $(\gS\cup\{1\})^n$ contains two elements $a,b$ which generates a free semigroup.
One of the main result proved in this note is:

\begin{thm}\label{thm1}
For any integer $d>1$ there is a constant $m=m(d)>1$ such that for any field $\BF$ and any finite subset $\gS\subset\GL_g(\BF)$ which generates a non-virtually solvable group we have $d^+(\gS)\le m$. It follows that there is a constant $C=C(d)>0$ such that $h(\gS)>C$ for any such $S$.
\end{thm}

The constants $m,C$ can be effectively estimate, but in order to keep this note short and simple we will allow ourselves to use some compactness arguments which will make the proofs non effective. Theorem \ref{thm1} answers a question of Gromov.
It improves the result of \cite{EMO} in the following three senses:
\begin{enumerate}
\item The analog result in \cite{EMO} was proved only for symmetric sets.
\item In \cite{EMO} it is assumed that the characteristic of $\BF$ is $0$.
\item The constants $m,C$ in Theorem \ref{thm1} do not depend on the actual group generated by $\gS$, but only on the dimension $d$.
\end{enumerate}

\begin{rem}
Uniform exponential growth for solvable non-virtually nilpotent groups was proved independently by Osin and by Wilson. The existence of a uniform constant $C(d)$ as in Theorem \ref{thm1} for solvable groups is very striking since it implies the well known Lehmer conjecture concerning the Mahler measure of algebraic integers. However for virtually solvable non-virtually nilpotent groups we only obtained the analog uniform result for discrete subgroups of $\GL_d(k)$ were $k$ is a local field. The details of this result will appear elsewhere.
\end{rem}

Using first order logic, it is straightforward to show that Theorem \ref{thm1} is equivalent to the following result concerning finite groups:

\begin{thm}
Given $d>1$, there is a function $f:\BN\to\BN$ tending to infinity, such that for any finite simple group $G$ of Lie type of dimension $d$, and any generating set $\gS$ of $G$ containing the identity, there are two elements $a,b\in\gS^{m(d)}$ such that the directed Cayley graph of $G$ with respect to $\{a,b\}$ has girth $\ge f(|G|)$.
\end{thm}

%------------------------------------------------------------------------------

\section{Eigenvalues vs. Norms}

Let $k$ be a local field with absolute value $|\cdot |_{k}.$ It induces the
standard norm on $k^{d}$ which in turn gives rise to an operator norm $%
\left\| \cdot \right\| _{k}$ on $M_{d}(k).$ If $k$ is not Archimedean, let $%
\mathcal{O}_{k}$ be its ring of integers and $m_{k}$ the maximal ideal in $%
\mathcal{O}_{k}.$ We note that $\left\| a\right\| _{k}\geq 1$ for all $a\in
\text{SL}_{d}(k).$ Let $\Lambda _{k}(a)$ be the maximum absolute value of
all eigenvalues of $a$ (recall that the absolute value has a unique
extension to the algebraic closure of $k$).

For a compact subset $Q\subset M_{d}(k)$ we denote:
\begin{eqnarray*}
\Lambda _{k}(Q) &=&\max \{\Lambda _{k}(a):a\in Q\} \\
\Vert Q\Vert _{k} &=&\max \{\Vert a\Vert _{k}:a\in Q\}~ \\
\Delta_{k}(Q) &=&\inf_{g\in \text{GL}_{d}(k)}\Vert gQg^{-1}\Vert _{k}.
\end{eqnarray*}
We will make use of the following strengthening of Lemma 4.2 from \cite{uti}:

\begin{lemma}
\label{lem:comparison}
There exists a constant $c=c(d)$ such that for any
compact subset $Q\subset M_{d}(k)$ we have
\begin{equation*}
\Lambda _{k}(Q^{i})\geq c\cdot \gD_k(Q)^{i}
\end{equation*}
for some $i\leq d^{2}$. Moreover, if $k$ is non-Archimedean the inequality
holds with $c=\left( |\pi |_{k}\right) ^{2d-1}$ for a uniformizer $\pi $ for
$k$.
\end{lemma}
\medskip

The proof of Lemma \ref{lem:comparison} relies on the following two lemmas. We refer to \cite{uti} for proofs of Lemmas \ref{nilp:span} and \ref{Lambda&Delta=0} .

\begin{lem}
\label{nilp:span}Let $R$ be a field or a finite ring and let $\mathcal{A}%
\leq M_{d}(R)$ be a subring and $R$-submodule. Suppose that $\mathcal{A}$ is
spanned as an $R$-module by nilpotent matrices, then $\mathcal{A}$ is
nilpotent, i.e. $\mathcal{A}^{n}=\{0\}$ for some $n\geq 1.$
\end{lem}

\begin{lem}
\label{Lambda&Delta=0} For a compact subset $Q\subset M_{d}(k)$ the
following are equivalent:

$(i)$ $Q$ generates a nilpotent subalgebra.

$(ii)$ $\gD_k(Q)=0$.

$(iii)$ $\Lambda _{k}(Q^{i})=0,~\forall i\leq d^{2}$.
\end{lem}

\begin{proof}[Proof of Lemmea \ref{lem:comparison}]
Suppose by contradiction that there is a sequence of compact sets
$Q_1,Q_2,\ldots$ in $M_d(k)$ such that $\Lambda_k(Q_n^i)<
{\gD_k(Q_n^i)}/{n},~\forall i\leq d^2$. By replacing $Q_n$ with a
suitable conjugate of it, we may assume that $\| Q_n\|_k\leq 2\gD_k
(Q_n)$, and by normalizing we may assume that $\| Q_n\|_k=1$. Let
$Q$ be a limit of $Q_n$ with respect to the Hausdorff topology on
$M_d(k)$. Then $\| Q\|_k=1$, $\gD_k (Q)\geq \frac{1}{2}$ since
$\gD_k$ is upper semi-continuous, and by continuity of
$\Lambda_k$, $\Lambda_k(Q^i)=0,~\forall i\leq d^2$. This however
contradicts Lemma \ref{Lambda&Delta=0}.

Let us now explain why one can take $c=\left( |\pi |_{k}\right)
^{2d-1}$ in the non archimedean case. Let $\mathcal{O}_k$ be the
ring of integers of $k$ and $m_k$ its maximal ideal. Suppose by
contradiction that there is a compact subset $Q$ of $M_{d}(k)$
such that $\Lambda _{k}(Q^{i})<c\cdot \gD_k(Q)^{i}$
for every $i\leq d^{2}$. Then up to conjugating by a suitable element of SL$%
_{d}(k)$ and renormalizing, we may assume that $\gD_k(Q)=\left\| Q\right\| =1
$. Note that since $\Lambda _{k}(Q^{i})\leq \left( |\pi |_{k}\right) ^{2d}$,
the coefficients of the characteristic polynomial of any element of $Q^{i}$
belong to $m_{k}^{2d}.$ Let $\mathcal{A}$ be the ring and $\mathcal{O}_{k}$%
-submodule generated by $Q$ in $M_{d}(\mathcal{O}_{k}).$ Let $\overline{%
\mathcal{A}}$ its image modulo $m_{k}^{2d}$ in $M_{d}(\mathcal{O}%
_{k}/m_{k}^{2d}).$ Then $\overline{\mathcal{A}}$ is spanned as an $\mathcal{O%
}_{k}/m_{k}^{2d}$-module by nilpotent matrices. Hence, according to Lemma
\ref{nilp:span}, $\overline{\mathcal{A}}$ is nilpotent and therefore lies in
the upper triangular matrices after a suitable change of basis. Hence $%
\mathcal{A}\leq
m_{k}^{2d}M_{d}(\mathcal{O}_{k})+N_{d}(\mathcal{O}_{k})$, where
$N_{d}(\mathcal{O}_{k})$ are the upper triangular matrices with
$0$ diagonal and coefficients in $\mathcal{O}_{k}.$ However, we
can now conjugate $\mathcal{A}$ by the matrix ${t}=diag(\pi
^{-(d-1)},\pi ^{-d-1},...,\pi ^{d-1})\in SL_{d}(k)$ and get that
${t}\mathcal{A}
{t}^{-1}\leq m_{k}^{2}M_{d}(\mathcal{O}_{k}),$ contradicting that $%
\gD_k(Q)=1.$%
%TCIMACRO{\TeXButton{End Proof}{\endproof}}
%BeginExpansion
\end{proof}

\begin{rem}\label{rem:c=1}
Note that in the non archimedean case, if we had considered the algebraic
closure $\overline{k}$ of $k$ instead of $k,$ then the constant $c$ could
have been taken to be equal to $1.$ Indeed, if $k_{n}$ denotes the
compositum of all extensions of $k$ of degree at most $n,$ and $\pi _{n}$ a
uniformizer for $k_{n},$ then $|\pi _{n}|_{\overline{k}}$ tends to $1$ as $n$
tends to infinity.
\end{rem}

%------------------------------------------------

%----------------------------------------------------------------------

\section{Lower bound on the algebraic entropy}

In this section, we prove the following theorem.

\begin{thm}
\label{thm:Q,delta,m} Let $G$ be a Zariski connected simple group
over a local field $k$, considered via its Adjoint representation
as an irreducible subgroup of $\SL (\mathfrak{g})$, and let
$\delta >0$. Then there is a constant $m=m(\dim G,\delta,k)$ such
that for any compact subset $Q$ generating a Zariski dense
subgroup, which is assumed to be non-compact in the
non-Archimedean case, with $\gD_k(Q)\geq 1+\delta $, we have
$d^{+}(Q)\leq m$.
\end{thm}

We shall require a few lemmas.

\begin{lem}
\label{N(C,delta,d)} For any $C$ and $\delta >0$ there is
$N(C,\delta ,d)$ such that if $Q\subset \text{SL}_{d}(k)$ is a
compact set with $\gD_k(Q)\geq 1+\delta $ then $\Vert \cup
_{i=1}^{N}Q^{i}\Vert _{k}>C$.
\end{lem}

\begin{proof}
Let us first assume that $k$ is Archimedean. Assuming the
contrary, we obtain a sequence of compact sets $Q_n$ with
$\gD_k(Q_n)\geq 1+\gd$ and $\|\cup_{i=1}^nQ_{n}^i\|\leq C$. Since
the norms of $Q_n$ are uniformly bounded we can pass to a
subsequence which converges to a limit compact set $Q$,
satisfying:
$$
 \|\cup_{i=1}^\infty Q^i\|\leq C.
$$
It follows that the closed semigroup generated by $Q$ is compact.
Since a compact sub-semigroup of a topological group is a group,
we obtain that $\gD_k(Q)=1$, contradicting the assumption that $\gD_k
(Q_n)\geq 1+\gd$, as $\gD_k$ is upper semi continuous.

Now if $k$ is non-Archimedean we can argue in the same way,
assuming further that $\langle Q_n\rangle$ are non-compacts, and
berrying in mind that any compact subgroup is contained in an open
compact subgroup.
\end{proof}

\begin{lem}\label{lem:proximal1}
Suppose $a\in SL_{d}(k)$ satisfies $\Lambda
_{k}(a)\geq 2$\label{conj}$\lambda _{k}(a)$ where $\lambda
_{k}(a)$ is the modulus of the second highest eigenvalue of $a$. Then the top eigenvalue $%
\alpha _{1}$ belongs to $k$, $|\alpha _{1}|=\Lambda _{k}(a)$ and
there exists $h\in SL_{d}(\overline{k})$ with $\left\| h\right\|
\leq 6^{d}\left\|
A\right\| ^{2d}$ such that the matrix $a^{\prime }=hah^{-1}$ is such that $%
a^{\prime }(e_{1})=\alpha _{1}e_{1}$ and $a^{\prime }(H)=H$ where $%
H=\left\langle e_{2},...,e_{d}\right\rangle $ and $\|a'|_{H}\|
\leq \frac{3}{2}\lambda _{k}(A).$
\end{lem}

Let us first prove:

\begin{sublem}
\label{lem:proximal2} For any $a\in SL_{d}(k)$ there exists $h\in
SL_{d}(\overline{k})$ with $\left\| hah^{-1}\right\| \leq
\frac{3}{2}\Lambda (a)$ and $\left\| h\right\| \leq 2^{d-1}\left\|
a\right\| ^{d-1}.$
\end{sublem}

\proof

We can triangularize $a$, i.e. find a basis $u_{1},...,u_{d}$ of $\overline{k%
}^{d}$ such that $a(u_{i})=\lambda _{i}u_{i}+$ other terms in $%
span\{u_{j},j>i\}.$ Then we can ``orthonormalize'' this basis,
i.e. we can find $g$ in the standard maximal compact subgroup of
$SL_{d}(\overline{k})$ (i.e. $\left\| g\right\| =1$) such that
$u_{i}=\mu _{i}ge_{i}+$ other terms
in $span\{u_{j},j>i\}$, where $\mu _{i}\in \overline{k}$ and $%
e_{1},...,e_{d} $ is the canonical basis of $\overline{k}.$ Hence
up to conjugating $a$ by $g^{-1}$ (this doesn't change the norm)
and rescaling each $ u_{i}$ we may assume that we had started with
$a$ triangular with coefficients in $\overline{k}.$

Now we can take $h=t^{\frac{d+1}{2}}diag(t^{-1},...,t^{-d})\in SL_{d}(%
\overline{k}).$ Let $a=(a_{ij})$, then $\max |a_{ij}|\leq \left\|
a\right\| . $ We get $\left\| hah^{-1}\right\| \leq \Lambda
(a)+2t^{-1}\left\| a\right\| .$ Indeed for $x\in \overline{k}$ we
have $hah^{-1}x=\sum \lambda
_{i}x_{i}e_{i}+\sum_{j>i}x_{i}t^{i-j}a_{ij}e_{j}.$ Let us estimate
the second term (say in the archimedean case):

\begin{eqnarray*}
\| \sum_{j>i}x_{i}t^{i-j}a_{ij}e_{j}\| ^{2} =\sum_{1\leq
i,i^{\prime }\leq d-1}x_{i}\overline{x_{i^{\prime
}}}\sum_{j>i,i^{\prime }}t^{i+i^{\prime
}-2j}a_{ij}\overline{a_{i^{\prime }j}} \leq \| a\| ^{2}\sum_{1\leq
i,i^{\prime }\leq d-1}|x_{i}||x_{i^{\prime }}|\sum_{j>i,i^{\prime
}}t^{i+i^{\prime }-2j}.
\end{eqnarray*}
However $\sum_{j>i,i^{\prime }}t^{i+i^{\prime }-2j}\leq t^{i+i^{\prime }}%
\frac{t^{-2(\max \{i,j\}+1)}}{1-t^{-2}}\leq 2t^{-2}t^{-|i-i^{\prime }|}$ if $%
1-t^{-2}\geq 1/2.$ So
\begin{equation*}
\| \sum_{j>i}x_{i}t^{i-j}a_{ij}e_{j}\| ^{2}\leq 2t^{-2}\| A\|
^{2}\sum_{1\leq i,i^{\prime }\leq d-1}|x_{i}||x_{i^{\prime
}}|t^{-|i-i^{\prime }|}.
\end{equation*}
The right hand side is the quadratic form associated to the matrix $%
B:=(t^{-|i-i^{\prime }|})_{i,i^{\prime }}.$ We compute easily
$\left\| B\right\| \leq 1+\frac{2t^{-1}}{1-t^{-1}}\leq 2$ if
$t\geq 3.$ Hence
\begin{equation*}
\| \sum_{j>i}x_{i}t^{i-j}a_{ij}e_{j}\| ^{2}\leq 4t^{-2}\| A\|
^{2}\| x\| ^{2}.
\end{equation*}

Thus we can take $t=4\left\| a\right\| /\Lambda (a)$ in case $k$
is archimedean. When $k$ is non-archimedean, we can take $t=\pi ^{-n},$ where $%
\left| \Lambda (a)\right| =\left| \pi ^{n-1}\left\| a\right\|
\right| $ (so that $\left| t\right| \leq \left\| a\right\|
^{2}/\Lambda (a)^{2}$). Hence in both cases we get $\left\|
hah^{-1}\right\| \leq \frac{3}{2}\Lambda (a)$ and $\left\|
h\right\| \leq t^{\frac{d-1}{2}}\leq 2^{d-1}\left\| a\right\|
^{d-1},$ since $\Lambda (a)\geq 1$ because $a\in
SL_{d}(\overline{k}).$
\endproof

\begin{proof}[Proof of Lemma \ref{lem:proximal1}] As in the proof of
sublemma \ref {lem:proximal2}, we can conjugate $a$ to some lower
triangular matrix without changing the norm of $a$. Let $\alpha
_{1},\alpha _{2},...,\alpha _{d}$ be the eigenvalues of $a$ in
$\overline{k},$ where we chose $|\alpha _{1}|=\Lambda _{k}(a)$ and
$|\alpha _{2}|=\lambda _{k}(a)$. Note that since $\Lambda
_{k}(a)>\lambda _{k}(a)$ all Galois conjugates of $\alpha _{1}$
(whose modulus is the same as that of $\alpha _{1}$) must be equal to $%
\alpha _{1}$ itself. Hence $\alpha _{1}\in k.$ Let $H$ be the hyperplane in $%
\overline{k}^{d}$ spanned by $e_{2},...,e_{d}$ and let $%
v=(1,x_{2},...,x_{d}) $ be the eigenvector corresponding to
$\alpha _{1}.$ We claim that $d(v,H)\geq 1/(1+L)^{d-1},$ where
$L=2\left\| a\right\| /\Lambda _{k}(a)$.

Indeed we have $d(v,H)\geq 1/\sqrt{1+\sum |x_{i}|^{2}}$ with
equality in the archimedean case. But if $(\alpha
_{1},a_{2},...,a_{d})$ is the first column of $a$, then
$a_{2}+\alpha _{2}x_{2}=\alpha _{1}x_{2}$ hence $|x_{2}|\leq
|a_{2}|/|\alpha _{1}-\alpha _{2}|\leq 2\left\| a\right\| /\Lambda
_{k}(a)$.
Then similarly we get $|x_{3}|\leq L+L^{2}$ and finally by induction $%
|x_{k}|\leq L(1+L)^{k-2}.$ We compute $1+\sum |x_{i}|^{2}\leq
(1+L)^{2d-2},$ hence $d(v,H)\geq 1/(1+L)^{d-1}.$

Now let $h_{1}\in SL_{d}(\overline{k})$ be defined by $h_{1}v=e_{1}$ and $%
h_{1}e_{i}=e_{i}$ for $i>1.$ It is easy to check directly in both
archimedean and non-archimedean case that $\left\| h_{1}\right\|
^{2}\leq
2(1+\sum |x_{i}|^{2}),$ hence $\left\| h_{1}\right\| \leq \sqrt{2}%
(1+L)^{d-1}.$

Finally, we can apply Lemma \ref{lem:proximal2} to the restriction
of $a$ to $H$ to get $h_{0}\in SL_{d}(\overline{k})$ with $\|
h_{0}\| \leq 2^{d-2}\| a\| ^{d-2}$ and $\left\|
h_{0}a_{|H}h_{0}^{-1}\right\|
\leq \frac{3}{2}\lambda _{k}(a).$ Finally letting $h=h_{0}h_{1}$ and $%
a^{\prime }=hah^{-1}$, we have $a^{\prime }e_{1}=\alpha _{1}e_{1},$ $%
a^{\prime }H=H,$ $\left\| a'|_{H}\right\| \leq \frac{3}{2}\lambda
_{k}(a)$ and $\left\| h\right\| \leq \left\| h_{0}\right\| \left\|
h_{1}\right\| \leq 6^{d}\left\| a\right\| ^{2d}.$
\end{proof}

The following lemma is verified by a direct simple computation. We
refer the reader to \cite{BG} Section 3, for a detailed study of
the dynamics of projective transformations.

\begin{lem}\label{lem:prox1}
Let $a\in\SL_d(k)$ and assume that $a(e_1)=\ga e_1$, $a(H)=H$,
where $H$ is the hyperplane spanned by $e_2,\ldots,e_d$, and $\|
a|_H\|\leq\gep^2|\ga |$, for some $\gep\leq 1/4$. Then $[a]$ is
$\gep$ proximal with respect to $([e_1],[H])$, i.e. it takes the
complement of the $\gep$-neighborhood of $[H]$ into the
$\gep$-neighborhood of $[e_1]$. Moreover $a$ is $\gep$-Lipschitz
outside the $\gep$-neighborhood of $[H]$.
\end{lem}

The following is a variant of the classical ping-pong lemma.

\begin{lem}[The Ping lemma]\label{lem:ping}
Let $a$ be as in Lemma \ref{lem:prox1}, and let $b\in\SL_d(k)$ be
another element with $b(e_1)\neq e_1$, $(\|
b\|\cdot\|b^{-1}\|)^2\leq\gep^{-\frac{1}{2}}-1$ and $d(b\cdot
[e_1],[H])\geq\sqrt{\gep}$. Then $a$ and $ba$ generate a free
semigroups.
\end{lem}

\begin{proof}
Let ${B}_\gep([e_1])$ denote the open $\gep$ ball around $[e_1]$.
The conditions on $b$ implies that $b({B}_\gep([e_1]))$ does not
meet the $\gep$-neighborhood of $[H]$, indeed by Lemma 3.1 in
\cite{BG} $b$ acts on $\BP^{d-1}$ with Lipschitz constant $(\|
b\|\|b^{-1}\|)^2$. Let
$$
 \gt =\sup \{ t\leq\gep :b(B_t([e_1]))\cap B_t([e_1])=\emptyset\}.
$$
Then $\gt>0$ and if we put $U=B_\gt ([e_1])$ and $V=b\cdot U$,
then $U\cap V=\emptyset$, $a\cdot (U\cup V)\subset U$ and $ba\cdot
(U\cup V)\subset V$.
\end{proof}

We shall also make use of the following:

\begin{lem}\label{lem:good-conjugation}
Let $G$ and $Q$ be as in Theorem \ref{thm:Q,delta,m} and consider
an irreducible representation $\rho$ of $G$ on $k^{d}$. Then there
is an element $q\in Q^{d+1}$ such that $\rho (q)\cdot [e_1]\neq
[e_1]$ and $\rho (q)\cdot [e_1]\notin [H]$.
\end{lem}

\begin{proof}
As $G$ is Zariski connected $\rho$ is strongly irreducible on
$k^d$ and hence there are $d-1$ elements in $Q^{d-1}\cdot [e_1]$
which together with $e_1$ form a basis to $k^d$. Denote the
corresponding projective points by $[e_1],q_1\cdot
[e_1],\ldots,e_{d-1}\cdot [e_1]$ where $q_i\in Q^{d-1}$. Since $G$
is Zariski connected it has no non-trivial map to the symmetric
group on $d$ elements and hence for some $q_d\in Q^d$ the
cardinality of the set $\{ [e_1],q_1\cdot [e_1],\ldots,q_{d}\cdot
[e_1]\}$ is $d+1$. If for some $i\leq d$ the point $q_i\cdot
[e_1]$ is not in $[H]$ then we are done. Assume that this is not
the case, then one can easily verify that our points $q_i\cdot
[e_1],~i=1,\ldots,d$ are in general position in $[H]$, i.e. the
lines corresponding to any $d-1$ of them spans $H$. Let $q_0\in Q$
be any element such that $q_0\cdot [H]\neq [H]$. Then for some
$i\neq j\leq d$, both $(q_0q_i)\cdot [e_1]$ and $(q_0q_i)\cdot
[e_1]$ are not in $[H]$, and at least one of them is different
from $[e_1]$.
\end{proof}

Equipped with the lemmas above we can now prove Theorem
\ref{thm:Q,delta,m}:

\begin{proof}[Proof of Theorem \ref{thm:Q,delta,m}]
We assume $d=\dim G$ and view $G$ as an irreducible subgroup of
$\SL_d(k)$. We may assume that $1\in Q$ and hence $Q^{j}\supset Q$
for any integer $j$, and in particular $Q^{j}$ generates a Zariski
dense subgroup.

Let $c$ be the constant from Lemma \ref{lem:comparison},
let $N=N(2^{d-1}/c,\delta ,d)$ be the constant from Lemma
\ref{N(C,delta,d)}. Replacing $Q$ by $Q^{N}$ we may assume that
$\gD_k(Q)\geq 2^{d-1}/c$.

Fix $i\leq d^{2}$ as in Lemma \ref{lem:comparison}, so that
$\Lambda _{k}(Q^{i})\geq c\gD_k(Q^{i})$ and take $a\in Q^{i}$ with
$\Lambda _{k}(a)=\Lambda _{k}(Q^{i})$. Then
\begin{equation*}
\Lambda _{k}(a)\geq c\gD_k\gD_k(Q^{i})\geq 2^{d-1},
\end{equation*}
and up to taking a suitable conjugation we may also assume
$\Lambda _{k}(a)\geq \frac{c}{2}\Vert Q^{i}\Vert $.

Taking a suitable irreducible factor of some wedge representation
(of degree at most $d$) we may assume that the top eigenvalue
$\alpha _{1}$ of $a$ is unique and that $|\alpha _{1}/\alpha
_{2}|_{k}\geq \Lambda_k(a)^{\frac{1}{d-1} }\geq 2$ where $\alpha
_{2}$ is second highest eigenvalue of $a$, so we can use Lemma
\ref{lem:proximal1} and find a suitable conjugating element $h\in
\text{SL}_{D}(k)$, where $D$ is the dimension of the new
representation. Note that the norm of a group element might have
increased by at most a power $d$ when we took the wedge
representation, hence the inequality above is replaced by
\begin{equation*}
\Vert a\Vert _{k}\geq \Lambda_{k}(a)\geq (c\Vert Q^{i}\Vert )^{\frac{1}{d}%
},~\text{~additionaly~}~\Vert h\Vert \leq 6^{D^2}\left\| a\right\|
^{2D^{2}}.
\end{equation*}
After conjugating with $h$ we may assume that $a(e_{1})=\alpha
_{1}e_{1},~a(H^{-})=H^{-}$, where $H^{-}=\text{span}\{e_{j}:j\geq
2\}$ and
\begin{equation*}
\Lambda_k(a)^{d}\geq \frac{c}{6^{2d^2}}\Vert Q^{i}\Vert
^{\frac{1}{4D^{5}}}.
\end{equation*}
Let $l$ be such that $|\ga_1/\ga_2|_k^l\geq \| Q^{d^2}\|$.

Let us first explain the proof under the assumption that $Q$ is
symmetric, and then give an alternative argument which holds also
without this assumption. By lemma \ref{lem:good-conjugation} we
can find $q\in Q^d$ such that $q\cdot [e_1]\neq [e_1]$ and $q\cdot
[e_1]\notin [H]$. It follows that the element $b=q^{-1}aq$
satisfies $b^n\cdot [e_1]\to q^{-1}\cdot [e_1]$ and in particular,
no power of $b$ stabilizes $[e_1]$. By the Cayley-Hamilton
theorem, we have $ \sum_{1\leq j\leq D}\beta _{j}b^{j}=-1$ for
some coefficients $\beta _{j}$ such that $|\beta _{j}|_{k}\leq
2^{D}\Lambda _{k}(b)^{D}$. Hence for some $1\leq j\leq D$, the
first coordinate of $b^{j}e_{1}$ is at least $\frac{1}{
D2^{D}\Lambda _{k}(b)^{D}}$. It follows that
$d([b^{j}e_{1}],H^{-})\geq \frac{1}{D2^{D}\Lambda _{k}(b)^{D}\Vert
b^{j}\Vert }$. Since the norm of $\| b\|$ is bounded by
$|\ga_1/\ga_2|^l$ we can take abounded power of $a$ so that it
will be $\gep$-proximal for suitable $\gep$ so that Lemma
\ref{lem:prox1} and Lemma \ref{lem:ping} will hold for $a,b$.

Let us now explain how to carry the proof without the assumption
that $Q$ is symmetric. All we need is to show that there is an
element $b$ in some bounded power of $Q$ such that $b^i\cdot
[e_1]\neq[e_1],~\forall i\leq D$. The existence of such an element
can be deduced from Lemma \ref{lem:Bezout}.

%Finally, note that it follows from the argument above the constant
%$m$ from the statement of \ref{thm:Q,delta,m} depends only on
%$\gd,\dim (G)$ and the local field $k$. We shall now show that it
%is actually independent of $k$. Indeed, the only places that there
%was a dependent on $k$ in the proof above is when we used
%Proposition \ref{prop:constant-c} and Lemma \ref{N(C,delta,d)}.
%However, as noted in Remark \ref{rem:c=1}, in the non Archimedean
%case we can assume that $c=1$, then it follows that
%$\Lambda_k(a)>1$ and since each eigenvalue of $a$ belongs to an
%extension of $k$ of degree at most $d$ it follows that
%$\Lambda_k(a)\geq |\pi |^{-\frac{1}{d}\geq$ where $\pi$ is an
%uniformaizer of $k$.
\end{proof}

The following lemma from \cite{uti} is a generalization of the corresponding lemma from \cite{EMO} and is proved by the same reasoning.

\begin{lem}\label{lem:Bezout}
Given an integer $\chi $ there is $N=N(\chi )$ such
that for any field $K$, any integer $d\geq 1$, any $K$--algebraic subvariety
$X$ in $GL_{d}(K)$ with $\chi (X)\leq \chi $ and any subset $\Sigma \subset {%
GL_{d}(K)}$ which contains the identity and generates a subgroup which is
not contained in $X(K)$, we have $\Sigma ^{N}\nsubseteq X(K)$.
\end{lem}

%--------------------------------------------------------------------

\section{Entropy for discrete subgroups}

In this section we prove the following:

\begin{thm}\label{thm:discrete}
Given an integer $d>1$, there is a constant $m'=m'(d)$ such that if $\gS\subset\GL_d(\BR)$ is a finite set, generating a discrete non-virtually solvable group, then $d^+(\gS)\le m'$.
\end{thm}

The proof relies on the classical Margulis lemma:

\begin{lem}\label{lem:Margulis-lem}
There exists a constant $\tau =\tau (d)>0$ such that if $\Sigma \subset
\text{SL}_{d}({\mathbb{R}})$ is a finite set generating a discrete non
virtually nilpotent group, then $E_{{\mathbb{R}}}(\Sigma )\geq 1+\tau $.
\end{lem}

The Margulis lemma is usually stated in terms of
displacement, however the exact same argument (c.f.
\cite{Thurston}) gives Lemma \ref{lem:Margulis-lem} with $\|\gS\|$
on the left side, and since discreteness is preserved under
conjugation, we may replace $\|~\|$ by $E_\BR$. Note also that the
analog statement hold for any local field, but, as remarked above,
we do not need it. We will assume in this section that $\Gamma
=\langle \Sigma \rangle $ is not virtually solvable. Let us first
reduce to the case of simple Zariski closure.

\begin{lem}
Let $\Gamma $ be a discrete non-virtually solvable linear group over $k$.
Then $\Gamma $ admits a representation $\rho $ with discrete image whose
Zariski closure is semisimple without compact factors.
\end{lem}

\begin{proof}
Let $G=\overline{\gC}^Z$ be the Zariski closure of $\gC$ and let
$A$ be the amenable radical of $G$, i.e. the maximal closed normal
amenable subgroup. Then $G/A$ is semisimple without compact
factors. Let $\rho$ be the restriction of the quotient map
$\rho:\gC\to G/A$. Since $\gC$ is discrete in $G$, the action of
$\gC$ on $G$ by left multiplications is amenable, and since $A$ is
amenable the action of $\gC$ on $G/A$ is also amenable.

In case $k=\BR$ it follows from Zimmer's theorem (c.f. \cite{BG1})
that the identity component of $\rho (\gC )$ is solvable, and
since it is also normal, as $\rho (\gC )$ is Zariski dense, it
must be trivial by the maximality of $A$, hence $\rho (\gC )$ is
discrete.

Assume now that $k$ is non-Archimedean and, by a way of
contradiction, assume that $\rho (\gC )$ is not discrete. By the
non-connected version of Zimmer's theorem (the generalized
Connes-Sullivan conjecture, c.f. \cite{BG1}), $\rho (\gC )$
contains an open solvable subgroup $\gC_0$. Let $U_0$ be its
closure, and let $U_n\subset U_0$ be a decreasing sequence of open
compact subgroups in $G/A$ with $\cap U_n=1$. Our assumption
implies that $U_n\cap\gC$ is non-trivial for every $n$. By
Noetherity of the Zariski topology the sequence
$\overline{U_n\cap\gC}^Z$ must stabilize after finitely many
steps, say $m$. It follows that $\overline{U_m\cap\gC}^Z$ is a
non-trivial closed normal solvable subgroup, contradicting the
maximality of $A$.
\end{proof}

\begin{proof}[Proof of Theorem \ref{thm:discrete}]
Let us assume now that $G=\overline{\Gamma }^{Z}$ is semisimple. Replacing $%
\Gamma $ by a subgroup of bounded index we may also assume that it is
Zariski connected. Indeed, look at the Adjoint representation of $G$. The
image of $\Gamma $ is discrete since the kernel is the center, hence
amenable. Moreover, the image of $\text{Ad}(G)$ is a subgroup of $\text{Aut}(%
\text{Ad}(G)^{\circ })$ and $|[\text{Aut}(\text{Ad}(G)^{\circ }):\text{Ad}%
(G)^{\circ }]|$ is bounded in terms of the number of factors of $\text{Ad}%
(G)^{\circ }$. Projecting to an appropriate factor of $G$ we can
also assume that it is simple. Hence we can apply Theorem
\ref{thm:Q,delta,m}.
\end{proof}

%-----------------------------------------------------------------------------

\section{Entropy of unbounded groups over non-Archimedean fields}

In this section we prove the following theorem:

\begin{thm}\label{thm:non-archimedean}
For any integer $d>1$ there is a number $m''=m''(d)$ such that for any non-archemedean local field $k$ and any subset $\gO\subset\GL_d(k)$ generating an unbounded group whose Zariski closure is semisimple, we have $d^+(\gO)\le m''$.
\end{thm}

The proof is similar to that of Theorem \ref{thm:Q,delta,m}, while instead of using the uniform gap $\gd$ we use the ultrametric inequality.

Suppose that $\gD_k(\gO)=\gb_k$. Since $\langle\gO\rangle$ is non-compact, $\gb_k>1$.
By Lemma \ref{lem:comparison} we may replace $k$ by some large algebraic extension $k'$ of $k$ so that the constant $c$ is arbitrarily close to $1$. We have to make sure however that $\gD_k(\gO)$ does not get arbitrarily close to $1$ in this procedure. This follows from the fact that the affine building $X_k$ associated to $\GL_d(k)$ embeds as a convex subset of the affine building $X_{k'}$ of $\GL_d(k')$, and as these spaces are CAT$(0)$ the projection from $X_{k'}$ to $X_k$ is $1$-Lipschitz. Since $\gO$, as a subset of $\SL_d(k)$, preserves $X_k$, it follows that its minimal displacement $\min(\gO):=\min\{ d(x,g\cdot x):g\in\gO,x\in X_{k'}\}$ is attained in $X_k$. Since $\min(\gO)$ and $\gD(\gO)$ are related by the following inequality (cf. \cite{uti}, Lemma 4.5):
$$
 \gD_k(\gO)\leq\exp(\min(\gO))\leq\gD_k^{\sqrt{d}},
$$
we see that $\gD_{k'}(\gO)$ is bounded away from $1$. Thus, replacing $k$ by an appropriate large extension, we may assume that the constant $c$ of Lemma \ref{lem:comparison} satisfies
$c^{-2}<\gb:=\gb_{k}$.

Assume that $1\in\gO$.
Let us also assume as we may by replacing $\gO$ by a subset of a bounded power of it, using Lemma \ref{lem:Bezout}, that all the elements in $\gO$ are semisimple. As in the proof of
Theorem \ref{thm:Q,delta,m} we may replace
the representation by a suitable irreducible factor of an exterior product, and obtain a strongly irreducible representation of dimension $D\le d^d$ and an element $a\in\gO^{d^2}$ with $|\ga_1/\ga_2|\ge\gb^\frac{1}{d-1}$ where $\ga_1,\ga_2$ are the eigenvalues of $a$ of highest and second highest absolute values.  Replace $\gO$ by a suitable conjugate, we may also assume that $\|(\gO)\|\le \gb^d$.

\begin{lem}
Let $v$ be the normalized eigenvector corresponding to $\ga_1$ and $H$ the hyperplane spanned by the other eigenvectors. Then $d([v],[H])\ge\frac{1}{\gb^{d^2}}$.
\end{lem}

\begin{proof}
Denote by $\OO_k$ the ring of integers in $k$.
By the Gram-Shmidt argument we can conjugate $a$ by an element of $\GL_D(\OO_k)$ and
obtain a triangular matrix whose first entry is $a_1$. Conjugating it further by a suitable diagonal matrix of norm $\le\gb^{d^2}$ we get an upper triangular matrix $a^g=gag^{-1}$ for which all the entries other the the first are of absolute value $<|a_1|$. This implies that $gH=\text{span}\{e_i:i\ge 2\}$, and hence the norm of the restriction $a^g|_{gH}$ is strictly less than $|\ga_1|$. This implies that $d([gv],[gH])=1$. Indeed, let $y\in gH$ be a normalized vector closest to $\hat{gv}$, and let $x=\hat{gv}-y$, then
$$
 |a_1|\cdot\| x\|=\| a^g\|\cdot\| x\|\ge\|a^g(x)\|=
 \| a^g(\hat{gv})-a^g(y)\|=\| a^g(\hat{gv})\|=|a_1|.
$$
The lemma follows.
\end{proof}

It follows that there is a matrix $h\in\GL_D(k)$ of norm $1$ with inverse of norm $\le \gb^{d^2}$ such that $h(e_1)=v$ and $h(\text{span}\{e_i:I>1\})=H$.

As in the proof of Theorem \ref{thm:Q,delta,m}, using the Cayley-Hamilton theorem, we find $b\in\gO^D$ such that $h(a^{d^{10}})h^{-1}$ and $h(ba^{d^{10}})h^{-1}$ play "ping" with open attracting sets $B(\gb^{d^5},[e_1])$ and $b^hB(\gb^{d^5},[e_1])$.
\qed

%-----------------------------------------------------------------------------

\section{Uniform entropy for linear groups over arbitrary fields}

We shall now give the proof of Theorem \ref{thm1}. By a global field
here we mean a finite algebraic extension of either $\Bbb{Q}$ or $\Bbb{F}%
_{q}(t),$ where $\Bbb{F}_{q}$ is the finite field with $q$ elements and $t$
is an indeterminate. Given a global field ${\mathbb{K}}$ and a finite set $S$
of places of ${\mathbb{K}}$ including all the infinite ones, we denote by ${%
\mathcal{O}}_{{\mathbb{K}}}(S)$ the ring of $S$--integers in ${\mathbb{K}}$.
The following lemma allows us to reduce the general case to Zariski dense subgroups of arithmetic groups.

\begin{lem}[cf. \cite{uti}]\label{lem:spec}
Let $\BF$ be a field and let $\Gamma\leq\GL_d(\BF)$ be a finitely generated group which is not virtually solvable. Then there is a subgroup $\gC'\leq\gC$ of index $\le d^2$, a global field ${\mathbb{K}}$, a finite set of places $S$ of ${\mathbb{K}}$ and a representation $f:\Gamma
^{\prime }\to \text{GL}_{d}({\mathcal{O}}_{{\mathbb{K}}}(S))$ whose image is Zariski dense
in a simple ${\mathbb{K}}$--algebraic group.
\end{lem}

Moreover, it is easy to check that if $\gS$ is a generating set for $\gC$ then $\Sigma ^{2i+1}$, where $i=[\gC:\gC']$, contains a generating set for $\gC'$.

In August 2006, Breuillard \cite{Br} announced the following beautiful Adelic version of the Margulis lemma:

\begin{thm}\label{thm:AML}
For every integer $d>1$ there is a constant $\gd=\gd_d>0$ such that for every number field $\BK$ and a finite set $\gS\subset\GL_d(\BK)$, which generates a group whose Zariski closure is simple, then either
\begin{itemize}
\item for some non-Archimedean completion $k$ of $\BK$ the group $\langle\gS\rangle$ is unbounded, or
\item for some Archimedean completion $k$ of $\BK$, $\gD_k(\gS)\ge 1+\gd$.
\end{itemize}
\end{thm}.

\begin{proof}[Proof of Theorem \ref{thm1}]
Let $\BF$ be a field and $\gS\subset\GL_d(\BF)$ a finite set generating a non-virtually solvable group. By Lemma \ref{lem:spec} it is enough to consider the case where $\gS\subset\GL_d(\BK)$ for some global field $\BK$ and $\overline{\langle\gS\rangle}^z$ is connected and simple. Thus by Theorem \ref{thm:AML} there is either a non-Archimedean completion $k$ of $\BK$ where $\langle\gS\rangle$ is unbounded, or an Archimedean completion $k$ of $\BK$ where $\gD_k(\gS)\ge 1+\gd(d)$. Thus Theorem \ref{thm1} follows from Theorems \ref{thm:Q,delta,m} and
\ref{thm:non-archimedean}. Note than in characteristic $p>0$ the proof does not rely on Theorem \ref{thm:AML}.
\end{proof}

%-----------------------------------------------------------------------------

\section{Uniform uniform Tits-alternative}

The uniform version of Tits alternative proved in \cite{uti} state that for every given finitely generated non-virtually solvable subgroup
$\gC\leq\GL_d(\BF)$ there is a constant $m_\gC$ such that for any symmetric subset $\gS$ of $\gC$ containing the identity which generates a Zariski dense subgroup (in particular every generating set), the $m_\gC$ ball $\gS^{m_\gC}$ contains two independent elements, i.e. two elements which generate a free non-abelian group.
Our next claim is that the constant $m_\gC$ could be made uniform for subgroups of $\GL_d$, i.e. that there is a uniform constant $m_d$ which applies for all non-virtually solvable subgroups of $\GL_d$. In order to avoid some technical difficulties and to keep this note short and simple we will give the proof only for $d=2$. The proof for general $d$ will be given in \cite{UUT}.

\begin{thm}\label{thm2}
There is a constant $m_2$ such that for any field $\BF$ and any
finite symmetric set $\gS\subset\SL_2(\BF)$ containing the
identity which generate a Zariski dense subgroup, the set
$\gS^{m_2}$ contains two independent elements.
\end{thm}

\begin{proof}
By Lemma \ref{lem:spec} we may replace $\BF$ by some global field $\BK$.

Assume first that $\text{char}(\BK)=0$. In that case if
$\BG<\BSL_2$ is a proper algebraic subgroup, then either
$\dim(\BG)=0$ and by Jordan's theorem $\BG$ admits an abelian
subgroup of bounded index, or $\dim(\BG)=1$ and $\BG^\circ$ is
unipotent and $\BG$ is solvable, or $\BG^\circ$ is diagonalizable
and and $[\BG:\BG^\circ]\le 2$, or $\dim(\BG)=2$ and $\BG$ is a
Borel subgroup, hence solvable. Thus a non Zariski dense subgroup
satisfies a simple law.

For any set $A$, denote by $A^{(2)}$ the set of squares $\{
a^2:a\in A\}$, and by $A'$ the set of commutators $\{ [a,b]:a,b\in
A\}$. By Lemma \ref{lem:comparison} if $n\ge 4$ then
$\gD(\gS^n)^2\ge\gD((\gS^n)^{(2)})\ge c\gD(\gS^n)^2$ for any
completion $k$ of $\BK$. By Lemma \ref{lem:Bezout}, up to
replacing $\gS$ by $\gS^{n_0}$ for some constant $n_0$, we may
assume that $\gS^{(2)\prime}$ (the set of commutators of squares)
does not satisfy the law from the previous paragraph, hence
generates a Zariski dense subgroup. We then, using Theorem
\ref{thm:AML}, chose an appropriate completion $k$ of $\BK$ which
is either non-Archimedean and $\langle\gS^{(2)\prime}\rangle$ is
unbounded, or Archimedean and $\gD_k(\gS^{(2)\prime})\geq 1+\gd$.
Let us assume that $k$ is Archimedean, the non-Archemedean case is
treated similarly.

Using Lemmas \ref{lem:comparison} and \ref{N(C,delta,d)} we get that for some constant $n_1$, $\gS^{n_1}$ contains a
semisimple element $a$ with eigenvalues $\ga,\frac{1}{\ga}$ with $|\ga|\ge 100$, such that up to conjugating $\gS$ by
a suitable element, using
Lemma \ref{lem:proximal1}, $a=\text{diag}(\ga,\frac{1}{\ga})$ and $\|\gS\|\le |\ga|^{f_1}$ for some constant $f_1$.
Set
$$
 d:=\min\{\max\{|\gs_{12}|:\gs\in\gS^{(2)}\},\max\{|\gs_{21}|:\gs\in\gS^{(2)}\}\}.
$$
The condition that $\|\gS^{(2)\prime}\|\ge 1+\gd$ implies that
$d\ge\ga^{-f_2}$ for some constant $f_2\ge f_1$. Indeed, otherwise
$\gS^{(2)'}$ would be "too" close to the set of upper (or lower)
triangular unipotents, and by conjugating it with suitable
diagonal element it would get too close to $1$.

Choose $b,c\in\gS$ such that $|(b^2)_{12}|,|(c^2)_{21}|\ge d$. Up to replacing $b,c$ by their inverses we may assume:
$$
 |\ga|^{-f_2}\le |b_{11}|,|b_{21}|,|c_{11}|,|c_{12}|\le |\ga|^{f_2}.
$$

Denote by $p=[(1,0)],q=[(0,1)]$ the attracting and repelling points of $a$. Since $d\ge|\ga|^{-f_2}$ it follows that for some $f_3$,
$$
 d(b\cdot p,\{p,q\})\ge|\ga|^{-f_3}~\text{~and~}~d(c^{-1}\cdot q,\{p,q\})\geq |\ga|^{-f_3}.
$$
Thus by multiplying a power of $a$ by $b$ on the left and by $c$ on the right, we get a very contracting element whose attracting and repelling points are bounded away from $p$ and $q$. In order to ensure that they are not close to each other we may multiply the new element further by a power of $a$ on both sides. This will "move" the attracting and repelling points of the new element close to $p$ and $q$ respectively, but with a distance which can be bounded from below. Hence if we chose a constant $n_5$ sufficiently large, and then a sufficiently larger constant $n_6$ and set
$$
 a'=a^{n_5}xa^{n_6}ya^{n_5},
$$
then $a,a'$ form a ping-pong pair, hence generate a free group.

\medskip

Suppose now that $\text{char}(\BK)>0$. The only difference is that we don't have an analog of Jordan's theorem for finite groups. Again using Lemma \ref{lem:Bezout} we may, after replacing $\gS$ by a bounded power of it, assume that $\gS^{(2)\prime}$ in addition to the condition above, does not generate a nilpotent subgroup. Then in case $\langle\gS^{(2)\prime}\rangle$ is infinite, it is Zariski dense and we can proceed as above. On the other hand if $\langle\gS^{(2)\prime}\rangle$ is finite then it is desecrate in any completion $k$ of $\BK$. We can then work in any $k$ in which $\langle\gS\rangle$ is unbounded since $\{\gS^{(2)\prime}\}$ cannot be "too" close to a unipotent subgroup, because if it was, the argument of Zassenhouse theorem would imply that it is unipotent. Thus we can apply the same proof as above also in this case. Note that in the positive characteristic the proof does not rely on Theorem \ref{thm:AML}.

\end{proof}

Again, by first order logic, Theorem \ref{thm2} could be formulated in terms of finite groups:

\begin{cor}
There is a function $f$ defined on prime powers $q=p^n$ and tending to infinity, such that for any symmetric generating set $\gS$ of $\SL_2(\BF_q)$ containing the identity, there are two elements $a,b\in\gS^{m_2}$ such that the Cayley graph $\chi(\SL_2(\BF_q),\{a,b\})$ has girth $\ge f(q)$.
\end{cor}

Theorem \ref{thm2} implies the following:

\begin{thm}\label{thm:UNA}
There is a constant $\gep_2>0$ such that if $A\subset\SL_2(\BF)$ is a set (not necessarily symmetric) generating a Zariski dense subgroup, and $B\subset\SL_2(\BF)$ is any finite set, then for some $a\in A$ we have
$$
 |aB\setminus B|\geq\gep_2|B|.
$$
\end{thm}

In terms of finite groups, Theorem \ref{thm:UNA} is formulated as:

\begin{cor}
Given $n,d$ there is $l=l(n,d)$ such that whenever $q=p^k\ge l$, for any $d$ generators $\gc_1,\ldots,\gc_d$ of $\SL_2(\BF_q)$ and any subset $A$ of $\SL_2(\BF_q)$ of size $\le n$ there is $1\le i\le d$ such that $\gc_i\cdot A\setminus A\ge\gep_2|A|$.
\end{cor}

%-------------------------------------------------------------------------------------------


\begin{thebibliography}{99}

\bibitem{Bo}  A. Borel, \textit{On free subgroups of semisimple groups}%
, L'Enseig. Math., t. \textbf{29} (1983) pp. 151--164.

\bibitem{Br}  E. Breuillard, preprint.

\bibitem{BG}  E. Breuillard, T. Gelander, \textit{On dense free subgroups of
Lie groups}, J. Algebra, \textbf{261} , no. 2, pp. 448--467, (2003).

\bibitem{BG1}  E. Breuillard, T. Gelander, \textit{A topological Tits
alternative}, to appear in Annals of Math.

\bibitem{note}  E. Breuillard, T. Gelander, \textit{Cheeger constant and
algebraic entropy of linear groups},  Int. Math. Res. Not.  2005,  no. 56, 3511--3523.

\bibitem{uti}  E. Breuillard, T. Gelander, \textit{Uniform independence in linear groups}, to appear in Invent. Math.

\bibitem{EMO}  A. Eskin, M. Mozes, H. Oh, \textit{On uniform exponential
growth for linear groups in characteristic zero}, Invent. Math. \textbf{160}%
, no. 1, pp. 1--30, (2005).

\bibitem{UUT} T. Gelander, \textit{A uniform uniform Tits alternative}, in preparation.

\bibitem{GH}  R. Grigorchuk, P. de la Harpe, \textit{Limit behaviour of
exponential growth rates for finitely generated group}s, in Essays on
geometry and related topics, Vol. 1, 2, 351--370, Monogr. Enseign. Math.,
\textbf{38}, (2001).

\bibitem{Thurston} W.P. Thurston, {\it Three-Dimensional Geometry and Topology},
Volume 1, Princeton univ. press, 1997.

\bibitem{Tits}  J. Tits, \textit{Free subgroups of Linear groups}, Journal
of Algebra \textbf{20}\textit{\ }(1972), 250-270.
\end{thebibliography}
\end{document}